\documentclass[11pt]{amsart}%
\usepackage{graphicx,amscd,color,amsmath,amsfonts,amssymb,geometry,amssymb}

\newtheorem{theorem}{Theorem}[section]

\begin{document}

\title[Grothendieck's theorem is optimal]{Grothendieck's theorem for absolutely summing multilinear operators is optimal}

\author[Pellegrino]{D. Pellegrino}
\address{Departamento de Matem\'{a}tica, \newline \indent
Universidade Federal da Para\'{i}ba, \newline \indent
58.051-900 - Jo\~{a}o Pessoa, Brazil.}
\email{dmpellegrino@gmail.com and pellegrino@pq.cnpq.br}

\author[Seoane]{J. B. Seoane-Sep\'{u}lveda}
\address{Departamento de An\'{a}lisis Matem\'{a}tico,\newline\indent Facultad de Ciencias Matem\'{a}ticas, \newline\indent Plaza de Ciencias 3, \newline\indent Universidad Complutense de Madrid,\newline\indent Madrid, 28040, Spain.}
\email{jseoane@mat.ucm.es}

\maketitle

\begin{abstract}
Grothendieck's theorem asserts that every continuous linear operator from
$\ell_{1}$ to $\ell_{2}$ is absolutely $\left(  1;1\right)  $-summing. In this
note we prove that the optimal constant $g_{m}$ so that every continuous
$m$-linear operator from $\ell_{1}\times\cdots\times\ell_{1}$ to $\ell_{2}$ is
absolutely $\left(  g_{m};1\right)  $-summing is $\frac{2}{m+1}$. We also show
that if $g_{m}<\frac{2}{m+1}$ there is $\mathfrak{c}$ dimensional linear space
composed by continuous non absolutely $\left(  g_{m};1\right)  $-summing
$m$-linear operators from $\ell_{1}\times\cdots\times\ell_{1}$ to $\ell_{2}.$
In particular, our result solves (in the positive) a conjecture posed by A.T.
Bernardino in 2011.
\end{abstract}

\section{Introduction}

A celebrated result of Grothendieck asserts that every continuous linear
operator from $\ell_{1}$ to $\ell_{2}$ is absolutely $\left(  1;1\right)
$-summing. It was recently proved \cite{ber} that this result can be lifted to
multilinear operators in the following fashion:

\begin{quote}
\textit{Every continuous $m$-linear operator from $\ell_{1}\times\cdots
\times\ell_{1}$ to $\ell_{2}$ is absolutely $\left(  \frac{2}{m+1};1\right)
$-summing.}
\end{quote}

In the same paper the author conjectured that the value $\frac{2}{m+1}$ is
optimal. A particular case of our main result gives a positive solution to
this conjecture:

\begin{theorem}
\label{yytt}The estimate $\frac{2}{m+1}$ is optimal. Moreover, if $g_{m}<\frac{2}{m+1}$ then there exists a $\mathfrak{c}$-dimensional linear space formed (except by the null vector) by continuous non absolutely $\left(g_{m};1\right)  $-summing $m$-linear operators. This result is optimal in terms of dimension.
\end{theorem}

Above, $\mathfrak{c}$ denotes the cardinality of the continuum. In other
words, our main result shows that if $g_{m}<\frac{2}{m+1}$, the set of
continuous non absolutely $\left(  g_{m};1\right)  $-summing multilinear
operators is $\mathfrak{c}$-lineable and moreover, maximal lineable. For the
theory of lineability we refer to \cite{ags,bams} and the references therein.

Our proof of the optimality of $\frac{2}{m+1}$ is inspired on ideas that date
back to the classical work of Lindenstrauss and Pe\l czy\'{n}ski \cite{lp}
and, later, explored in a series of papers (see, e.g., \cite{PAMS, BP,
mathZ,studia}).

Throughout this note, $X$ and $Y$ shall stand for Banach spaces over
$\mathbb{K}=\mathbb{R}$ or $\mathbb{C}$. The closed unit ball of $X$ is
denoted by $B_{X}$ and the topological dual of $X$ by $X^{\ast}$. Also, recall
that a continuous linear operator $u:X\rightarrow Y$ is absolutely
$(q;1)$-summing (see \cite{DJT}) if there exists $C\geq0$ such that
\[
\left(  \sum_{j=1}^{n}\Vert u(x_{j})\Vert^{q}\right)  ^{\frac{1}{q}}\leq
C\sup_{\varphi\in B_{X^{\ast}}}\sum_{j=1}^{n}|\varphi(x_{j})|
\]
for every $n\in\mathbb{N}$ and $x_{1},\ldots,x_{n}\in X$. The nonlinear theory
of absolutely summing operators was designed by Pietsch in 1983 (\cite{P}) and
since then has been intensively studied. One of the possible polynomial
generalizations of absolutely summing operators is the concept of absolutely
summing polynomial. The space of continuous $m$-homogeneous polynomials from
$X$ to $Y$ will be henceforth denoted by $\mathcal{P}(^{m}X;Y).$ Given a
positive integer $m$ and $q\geq\frac{1}{m}$, a continuous $m$-homogeneous
polynomial $P:X\rightarrow Y$ is absolutely $\left(  q;1\right)  $-summing if
there exists a constant $C\geq0$ such that
\[
\left(  \sum_{j=1}^{n}\Vert P(x_{j})\Vert^{q}\right)  ^{\frac{1}{q}}\leq
C\left(  \sup_{\varphi\in B_{X^{\ast}}}\sum_{j=1}^{n}|\varphi(x_{j})|\right)
^{m}%
\]
for every $n\in\mathbb{N}$ and $x_{1},\ldots,x_{n}\in X$. If $q<1,$ the
infimum of the constants $C$ satisfying the above inequality is a Banach
quasinorm for the space of absolutely $\left(  q,1\right)  $-summing
polynomials from $X$ to $Y$, and it is denoted by $\pi_{q,1}$. For multilinear
mappings the definition is similar:

A continuous $m$-linear operator $T:X\times\cdots\times X\rightarrow Y$ is
absolutely $\left(  q;1\right)  $-summing (with $q\geq\frac{1}{m}$) if there
is a constant $C\geq0$ such that
\[
\left(  \sum_{j=1}^{n}\Vert T(x_{j}^{(1)},...,x_{j}^{(m)})\Vert^{q}\right)
^{\frac{1}{q}}\leq C{\displaystyle\prod\limits_{k=1}^{m}}\left(  \sup
_{\varphi\in B_{X^{\ast}}}\sum_{j=1}^{n}|\varphi(x_{j}^{(k)})|\right)
\]
for every $n\in\mathbb{N}$ and $x_{1}^{(k)},\ldots,x_{n}^{(k)}\in X$, and
$k=1,...,m$. If $q<1,$ the infimum of the constants $C$ satisfying the above
inequality is a Banach quasinorm for the space of absolutely $\left(
q;1\right)  $-summing $m$-linear operators from $X\times\cdots\times X$ to
$Y.$

\section{The proof of Theorem \ref{yytt}}

Let $1\leq g_{m}<\frac{2}{m+1}$. The first part of our argument is mentioned
\textit{en passant}, without proof, in \cite{mathZ}, but since we have a more
self-contained approach, we present the details for the sake of completeness.
Let $n\in\mathbb{N}$ and $x_{1},\ldots,x_{n}\in\ell_{1}$ be non null vectors.
Consider $x_{1}^{\ast},\ldots,x_{n}^{\ast}\in B_{\ell_{\infty}}$ so that
$x_{j}^{\ast}(x_{j})=\left\Vert x_{j}\right\Vert $ for every $j=1,\ldots,n$.
Let $a_{1},\ldots,a_{n}$ be scalars such that $\sum\limits_{j=1}^{n}%
|a_{j}|^{\frac{2}{g_{m}}}=1$ and define
\[
P_{n}\colon\ell_{1}\longrightarrow\ell_{2}~,~P_{n}(x)=\sum\limits_{j=1}%
^{n}\left\vert a_{j}\right\vert ^{\frac{1}{g_{m}}}x_{j}^{\ast}(x)^{m}e_{j},
\]
where $e_{j}$ is the $j$-th canonical vector of $\ell_{2}.$ For every
$x\in\ell_{1}$,
\[
\left\Vert P_{n}(x)\right\Vert =\left(  \sum\limits_{j=1}^{n}\left\vert
\left\vert a_{j}\right\vert ^{\frac{1}{g_{m}}}x_{j}^{\ast}(x)^{m}\right\vert
^{2}\right)  ^{\frac{1}{2}}\leq\left(  \sum\limits_{j=1}^{n}\left\vert
a_{j}\right\vert ^{\frac{2}{g_{m}}}\right)  ^{\frac{1}{2}}\left\Vert
x\right\Vert ^{m}=\left\Vert x\right\Vert ^{m}.
\]

Since $P_{n}$ is a polynomial of finite type, then it is plain that $P_{n}$ is
absolutely $\left(  g_{m};1\right)  $-summing. Note that for $k=1,\ldots,n$,
we have
\[
\Vert P_{n}(x_{k})\Vert=\left\Vert \sum\limits_{j=1}^{n}\left\vert
a_{j}\right\vert ^{\frac{1}{g_{m}}}x_{j}^{\ast}(x_{k})^{m}e_{j}\right\Vert
\geq\left\vert a_{k}\right\vert ^{\frac{1}{g_{m}}}x_{k}^{\ast}(x_{k}%
)^{m}=\left\vert a_{k}\right\vert ^{\frac{1}{g_{m}}}\Vert x_{k}\Vert^{m}.
\]
To simplify the notation we write $\Vert(x_{j})_{j=1}^{n}\Vert_{w,1}%
:=\sup_{\varphi\in B_{X^{\ast}}}\sum_{j=1}^{n}|\varphi(x_{j}^{(k)})|$. Thus,
we have
\begin{align*}
\left(  \sum\limits_{j=1}^{n}\Vert x_{j}\Vert^{mg_{m}}\left\vert
a_{j}\right\vert \right)  ^{\frac{1}{g_{m}}} &  =\left(  \sum\limits_{j=1}%
^{n}\left(  \Vert x_{j}\Vert^{m}\left\vert a_{j}\right\vert ^{\frac{1}{g_{m}}%
}\right)  ^{g_{m}}\right)  ^{\frac{1}{g_{m}}}\\
&  \leq\left(  \sum\limits_{j=1}^{n}\Vert P_{n}(x_{j})\Vert^{g_{m}}\right)
^{\frac{1}{g_{m}}}\\
&  \leq\pi_{g_{m},1}(P_{n})\Vert(x_{j})_{j=1}^{n}\Vert_{w,1}^{m}.
\end{align*}
Since this last inequality holds whenever $\sum\limits_{j=1}^{n}|a_{j}%
|^{\frac{2}{g_{m}}}=1$, denoting $\left(  \frac{2}{g_{m}}\right)  ^{\ast}$ to the
conjugate of $\frac{2}{g_{m}}$ we obtain
\begin{align*}
\left(  \sum\limits_{j=1}^{n}\Vert x_{j}\Vert^{mg_{m}\left(  \frac{2}{g_{m}%
}\right)  ^{\ast}}\right)  ^{\frac{1}{\left(  \frac{2}{g_{m}}\right)  ^{\ast}%
}} &  \leq\sup\left\{  \sum\limits_{j=1}^{n}\left\vert a_{j}\right\vert \Vert
x_{j}\Vert^{mg_{m}};\sum\limits_{j=1}^{n}|a_{j}|^{\frac{2}{g_{m}}}=1\right\}
\\
&  \leq\left(  \pi_{g_{m},1}(P_{n})\Vert(x_{j})_{j=1}^{n}\Vert_{w,1}%
^{m}\right)  ^{g_{m}}%
\end{align*}
and, then,
\begin{equation}
\frac{\left(  \sum\limits_{j=1}^{n}\Vert x_{j}\Vert^{mg_{m}\left(  \frac
{2}{g_{m}}\right)  ^{\ast}}\right)  ^{\frac{1}{g_{m}\left(  \frac{2}{g_{m}%
}\right)  ^{\ast}}}}{\Vert(x_{j})_{j=1}^{n}\Vert_{w,1}^{m}}\leq\pi_{g_{m}%
,1}(P_{n}).\label{yt}%
\end{equation}
Since $1\leq g_{m}<\frac{2}{m+1}$ we have $mg_{m}\left(  \frac{2}{g_{m}%
}\right)  ^{\ast}<2$ and from a weak version of the \emph{Dvoretzky-Rogers
Theorem} we know that $id_{\ell_{1}}$ is not $\left(  mg_{m}\left(  \frac
{2}{g_{m}}\right)  ^{\ast};1\right)  $-summing. Combining this fact with
(\ref{yt}) we conclude that we can find $x_{j}$ in $\ell_{1}$ for all positive
integer $j$ so that
\begin{equation}
\lim_{n\rightarrow\infty}\pi_{g_{m},1}(P_{n})=\infty\text{ and }\left\Vert
P_{m}\right\Vert =1.\label{65}%
\end{equation}
We thus conclude that the space of all absolutely $\left(  g_{m};1\right)
$-summing $m$-homogeneous polynomials from $\ell_{1}$ to $\ell_{2}$ is not
closed in $\mathcal{P}\left(  ^{m}\ell_{1};\ell_{2}\right)  .$ In fact,
otherwise, since the quasinorm $\pi_{g_{m},1}$ is complete, the \emph{Open
Mapping Theorem} to $F$-spaces would contradict (\ref{65}).

Now, let $P:\ell_{1}\rightarrow\ell_{2}$ be a continuous non $\left(
g_{m};1\right)  $-summing $m$-homogeneous polynomial. Split $\mathbb{N}$ into
a countable union of pairwise disjoint countable sets $\mathbb{N}%
_{1},\mathbb{N}_{2},...$. For all $j$, let $$\mathbb{N}_{j}=\left\{
a_{1}^{(j)}<a_{2}^{(j)}<\cdots\right\},  $$ and define $P^{(j)}:\ell
_{1}\rightarrow\ell_{2}$ by $\left(  P^{(j)}(x)\right)  _{a_{k}^{(j)}}=\left(
P(x)\right)  _{k}$ and $\left(  P^{(j)}(x)\right)  _{k}=0$ if $k\notin
\mathbb{N}_{j}$. It is simple to prove that $P^{(j)}$ is also a continuous non
$\left(  g_{m};1\right)  $-summing $m$-homogeneous polynomial and the set
$\left\{  P^{(1)},P^{(2)},...\right\}  $ is linearly independent. Finally, we
note that the linear operator $\Phi:\ell_{1}\rightarrow\mathcal{P}\left(
^{n}\ell_{1};\ell_{2}\right)  $ given by $\left(  \beta_{j}\right)
_{j=1}^{\infty}\mapsto%
{\textstyle\sum\limits_{j=1}^{\infty}}
\beta_{j}P^{(j)}$ is injective and it is simple to verify that $\Phi\left(
\ell_{1}\right)  $ is composed (except by the null vector) exclusively by non
absolutely $\left(  g_{m};1\right)  $-summing $m$-homogeneous polynomials. We
thus conclude that the set of continuous non $\left(  g_{m};1\right)
$-summing $m$-homogeneous polynomials is $\Phi\left(  \ell_{1}\right)
$-lineable. Since $\dim\Phi\left(  \ell_{1}\right)  =\dim\ell_{1}%
=\mathfrak{c}$ we conclude that this set is $\mathfrak{c}$-lineable.

Since $\mathcal{P}\left(  ^{n}\ell_{1};\ell_{2}\right)  $ is isomorphic to the
space of symmetric $m$-linear operators from $\ell_{1}\times\cdots\times
\ell_{1}$ to $\ell_{2}$ and since $P$ is absolutely $\left(  g_{m};1\right)
$-summing if and only if its associated symmetric $m$-linear operator is
absolutely $\left(  g_{m};1\right)  $-summing, our result is translated to the
multilinear setting.

The above result is optimal in terms of dimension. In fact, it is well known
that $\ell_{1}$ is isometric to the completion of its projective tensor
product, i.e., $\ell_{1}=\ell_{1}\widehat{\otimes}_{\pi}\cdots\widehat
{\otimes}_{\pi}\ell_{1}.$ Thus%
\[
\dim\mathcal{L}\left(  ^{m}\ell_{1};\ell_{2}\right)  =\dim\mathcal{L}\left(
\ell_{1}\widehat{\otimes}_{\pi}\cdots\widehat{\otimes}_{\pi}\ell_{1};\ell
_{2}\right)  =\dim\mathcal{L}\left(  \ell_{1};\ell_{2}\right)  =\mathfrak{c.}%
\]

We remark that if $\pi_{g_{m},1}$ was locally convex, since we have proved
that the space of all absolutely $\left(  g_{m};1\right)  $-summing
$m$-homogeneous polynomials from $\ell_{1}$ to $\ell_{2}$ is not closed in the
space of all continuous $m$-homogeneous polynomials from $\ell_{1}$ to
$\ell_{2}$, then from a result due to Drewnowski (see \cite[Theorem 5.6 and
its reformulation ]{drew}) we would conclude that the set of all continuous
$m$-homogeneous polynomials from $\ell_{1}$ to $\ell_{2}$ that fail to be
absolutely $\left(  g_{m};1\right)  $-summing is spaceable, i.e., contains
(except for the null vector) a closed infinite-dimensional subspace.

\end{document}